\documentclass[11pt]{article}

\usepackage{graphicx}
\usepackage{amssymb, amsmath,}
                                                                                                                                                                                                                                                                                                                                                                                                                                                                                                                                                                                                                                                                                                                                                                                                                                                                                                                                                                                                                                                                                                                                               \newtheorem{theorem}{Theorem}[section]
\newtheorem{corollary}[theorem]{Corollary}
\newtheorem{definition}[theorem]{Definition}

\newtheorem{Remark}[theorem]{Remark}

\title{\textbf{New informations on the structure of the functional codes defined by forms of degree
 $h$ on non-degenerate Hermitian varieties in $\mathbb{P}^{n}(\mathbb{F}_q)$}}

\date{}
\author{$\text{Fr\'ed\'eric A. B. Edoukou}^{*}$, $\text{San Ling}^{*}$, Chaoping Xing{\thanks{This work is supported by
MOE-AcRF Tier 2 Research Grant, Singapore (No.T206B2204).}}\\
Division of Mathematical Sciences,\\  Nanyang Technological University, \\
  21 Nanyang Link,
 Singapore 637371. \\
 E.mail : $\{\mathrm{abfedoukou, lingsan, xingcp} \}$@ntu.edu.sg}
\begin{document}
\maketitle

{\footnotesize \begin{flushleft} \textbf{Abstract}
\end{flushleft}
We study the functional codes of order $h$ defined by G. Lachaud on
$\mathcal{X} \subset {\mathbb{P}}^n(\mathbb{F}_q)$ a non-degenerate
Hermitian variety. We give a condition of divisibility of the weights
of the codewords. For $\mathcal{X}$ a non-degenerate Hermitian surface,
we list the first five weights and the corresponding codewords and
give a positive answer on a conjecture formulated on this question.
The paper ends with a conjecture on the minimum distance  and the distribution
of the codewords of the first $2h+1$ weights of the  functional codes for the functional codes of
order $h$ on $\mathcal{X} \subset {\mathbb{P}}^n(\mathbb{F}_q)$ a non-singular Hermitian variety.
}
\section{Introduction}
Let $\mathcal{X}$ be a projective algebraic variety over the finite field
$\mathbb{F}_q$. The functional codes $C_h(\mathcal{X})$ defined by evaluating
the polynomials functions over the rational points of the algebraic variety
 $\mathcal{X}$ have been studied in general way by several authors. The works
 of Goppa on codes constructed on the non-singular Hermitian curves inspired  also
 several authors mainly I. M. Chakravarti and his group at the university
 of North-Carolina,  to generalize theses codes to the non-singular
 Hermitian surfaces. For those who first of all, want to have a very readable treatment of the works Goppa on
 Hermitian curves but have a limited understanding of algebraic geometry the book of
 W. Cary and V. Pless [2, pp.526-544] can be examined.\\
 Some interesting results have been obtained in the case of the non-singular Hermitian surfaces by
 I. M. Chakravarti's group [3]. Their works were mainly done on the fields $\mathbb{F}_4$
 of order four. Therefore computer programs have been used to find all the structure of the codes.
 In 1991, A. B. S\o rensen, in his Ph. D. Thesis [19, pp.7-9] recovered a part of
 their results by tools of algebraic geometry and finite geometry, and try also to study
 the codes $C_h(\mathcal{X})$ over the non-singular Hermitian surface by increasing the
 degree of the polynomial functions and the order of the fields. In 2007, in his Ph.D. Thesis [5] the first author of this
 paper continue the works of A. B. S\o rensen and solved part of  A. B. S\o rensen's conjecture formulated in the
  past years on Hermitian surface. He generalized the results to the code $C_2(\mathcal{X})$ constructed on
  the non-singular Hermitian solid (non-singular Hermitian varieties of dimension three) where he found a
   structure on the
 repartition of the first five weights, their frequency and a divisibility condition for all the weights of the
 code $C_2(\mathcal{X})$. He also stated two conjectures one on the minimum weight codewords of the
  codes $C_2(\mathcal{X})$ constructed
 on  $\mathcal{X}: x_0^{t+1}+x_1^{t+1}+x_2^{t+1}+x_3^{t+1}+x_4^{t+1}=0$ where $h \le t$, and the second on the
repartition of the first five weights of the codes  $C_2(\mathcal{X})$ defined
on the non-singular Hermitian variety
$\mathcal{X}:x_0^{t+1}+x_1^{t+1}+...+x_{n-1}^{t+1}+x_n^{t+1}=0$ in
$\mathbb{P}^n(\mathbb{F}_q)$. The second conjecture formulated by him has been solved recently by
A. Hallez and L. Storme under the condition that $n< O(t^2)$. Under this restrictive condition we
can also remark a divisibility condition of the first five
weights which was not mentioned explicitly in their paper [9, p.9]. \\
The purpose of this paper is to give some news informations  on the
 structure of the codes $C_h(\mathcal{X})$. The paper has been organized as
 follows. First we recall some generalities on the number of solutions of a family of polynomials over a finite field.
Secondly we give the information on the divisibility condition
which should be respect by all the weights of the codes $C_h(\mathcal{X})$
constructed on the non-singular Hermitian variety :
$\mathcal{X}:x_0^{t+1}+x_1^{t+1}+...+x_{n-1}^{t+1}+x_n^{t+1}=0$ in
$\mathbb{P}^n(\mathbb{F}_q)$ with $q=t^2$ ($t$ is a prime power)
even if, theses weights are not computed explicitly. Thirdly by
using this result on the divisibility condition, we also solved the
conjecture formulated in [5, p.55], [7, p.113] on the fourth and
fifth weights of the codes $C_h(\mathcal{X})$ defined on the non-singular
Hermitian surface $\mathcal{X}:
x_0^{t+1}+x_1^{t+1}+x_2^{t+1}+x_3^{t+1}=0$ in
$\mathbb{P}^3(\mathbb{F}_q)$ with $q=t^2$ ($t$ is a prime power).\\
The paper ends with a conjecture on the minimum distance and
the distribution of the first $2h+1$ weights of the code
$C_h(\mathcal{X})$ where $\mathcal{X} \subset
{\mathbb{P}}^n(\mathbb{F}_q)$  is the non-singular Hermitian
variety.

 \section{Generalities}
We denote by $\mathbb{F}_q$ the field with $q$ elements. Let
$V=A^{n+1}(\mathbb{F}_q)$ be the affine space of dimension $n+1$
over $\mathbb{F}_q$ and ${\mathbb{P}}^{n}(\mathbb{F}_q)=\Pi_n$ the
corresponding projective space. Then
$$\pi_n=\#{\mathbb{P}^{n}(\mathbb{F}_q)}=q^n+q^{n-1}+...+1.$$
We denote by
$W_i$ the set of points with homogeneous coordinates
$(x_0:...:x_n)\in \mathbb{P}^n(\mathbb{F}_q)$ such that $x_j=0$ for
$j<i$ and $x_i \ne 0$. The family $\{W_i\}_{0\le i\le n}$ is a
partition of $\mathbb{P}^n( \mathbb{F}_q)$.
We use the term forms of degree $h$ to describe homogeneous
polynomials $f$ of degree $h$, and $\mathcal{V}=Z(f)$ (the zeros of $f$ in the projective space
 $\mathbb{P}^{n}(\mathbb{F}_q)$) is a hypersurface of degree $h$.
Let $\mathcal{F}_{h}$ be the vector space of forms of degree $h$ in
$V$. For any polynomial $f \in \mathcal{F}_{h}$ and
any point $P \in \mathbb{P}^{n}(\mathbb{F}_q)$ we define
  $$f(P)= f(x_0,...,x_n)/{x_i}^h \qquad \mathrm{with}\qquad P=(x_0:...:x_n) \in W_i$$
\begin{theorem} \label{TSS}(Tsfasman-Serre-S\o rensen)(\lbrack 18, p.351\rbrack, \lbrack 19, chp.2, pp.7-10\rbrack)
Let $f(x_{0},...,x_{n})$ be a homogeneous polynomial in $n+1$ variables with coefficients in
$\mathbb{F}_{q}$  and degree $h\le q$. Then the number of zeros of $f$ in  $\mathbb{P}^{n}(\mathbb{F}_q)$
satisfies:  $$\#Z_{(f)}(\mathbb{F}_{q})\le hq^{n-1}+\pi_{n-2}.$$
This upper bound is attained when $Z_{(f)}$ is a union of $h$ hyperplanes passing through a common
linear space of codimension 2.
 \end{theorem}
\begin{theorem}\label{quadrique type}(Ax-Katz)[13, p.485] Let $\mathbb{F}_q$ be a finite fields
 of characteristic p, having $q=p^a$ elements. Let N(S,T,f) defined as the number
 of points of V(S,T,f) with values in $\mathbb{F}_q$ and $\lambda$(S,T,f) defined as the least non-negative
 integer which is greater than
 $$\frac{Card(S)-\sum{d_i}}{\sup{(d_i)}}.$$
 Then N(S,T,f) $\equiv$ 0 modulo $q^{\lambda(S,T,f)}$.
\end{theorem}
Let  $\mathcal{X}\subset \mathbb{P}^{n}(\overline{\mathbb{F}_q})$ an algebraic variety and
$\# \mathcal{X}(\mathbb{F}_q )$ the
number of rational points of $\mathcal{X}$ over $\mathbb{F}_q$. The code $C_h(\mathcal{X})$ is the
image of the linear map
 $c: \mathcal{F}_{h}
  \longrightarrow
  \mathbb{F}_{q}^{\# \mathcal{X}(\mathbb{F}_q )}$, defined by $c(f)={(c_x(f)})_{x\in X}$,  where
   $c_x(f)= f(x_0,...,x_n)/{x_i}^h$ with $x=(x_0:...:x_n) \in W_i$. The length of
   $C_h(\mathcal{X})$ is equal to $\#\mathcal{X}(\mathbb{F}_q)$. The dimension of $C_h(\mathcal{X})$ is equal to
    $\dim{{\mathcal{F}}_h}-\dim{\ker c}$. Therefore, when $c$ is injective we get:
  \begin{equation}
  \label{dimducode}
  \dim\ {C_{h}(\mathcal{X})} =\left(
  \begin{array}{c}
  n+h\\
  h
\end{array}
\right).
\end{equation}
The minimum distance of $C_h(\mathcal{X})$ is equal to the minimum over all
$f$ of $\#\mathcal{X}(\mathbb{F}_q)-\#{\mathcal{X}}_{ Z(f)}(\mathbb{F}_{q})$.
\section{Divisibility condition on the weights of the codes  $C_h(\mathcal{X})$ defined on the
non-degenerate Hermitian variety}
In this section $\mathbb{F}_q$ denotes the field with $q$ elements,
where $q=t^2$ and $\mathcal{X}$
 denotes the non-degenerate (i.e. non-singular) Hermitian variety of $\mathbb{P}^n(\mathbb{F}_q) $
  of equation $\mathcal{X}: x_0^{t+1}+x_1^{t+1}+...+x_{n-1}^{t+1}+x_n^{t+1}=0$.\\
In $\lbrack 1,\mathrm{p}.1175\rbrack$ R. C. Bose and I. M. Chakravarti proved
the following result:
\begin{theorem}\label{cardinal hermitienne}
Let $\tilde{\mathcal{X}}\subset \mathbb{P}^n(\mathbb{F}_{q})$ be a
non-degenerate Hermitian variety. Then,
\begin{equation}
\label{bochanbrehermitian}
\#\tilde{\mathcal{X}}(\mathbb{F}_{q})=\lbrack
t^{n+1}-(-1)^{n+1}\rbrack\lbrack t^{n}-(-1)^{n}\rbrack/(t^2-1)
\end{equation}
\end{theorem}

\subsection{Transformation of a non-degenerate Hermitian variety to a quadric variety}

Here we will give an important result on the code defined on the
Hermitian variety  in the $\mathbb{P}^n(\mathbb{F}_q) $. Let $\mathcal{X}
:x_0^{t+1}+x_1^{t+1}+...+x_{n-1}^{t+1}+x_n^{t+1}=0$ be the
non-degenerate Hermitian variety over the field
$\mathbb{F}_{t^2}$.\\
If $\phi$ is an Hermitian form in  $n+1$ variables with $n\ge3$ over
the field $\mathbb{F}_{q}=\mathbb{F}_{t^2}$ defining the Hermitian variety $\mathcal{X}$, we have:
$$\phi(x_0, ...,x_n)=x_0^{t+1}+...+x_n^{t+1}$$
Let us denote by  $\alpha$ an element of $\mathbb{F}_{q}$ which is
not in $\mathbb{F}_{t}$. One has $$\mathbb{F}_{q}=
\mathbb{F}_{t}\oplus\alpha \mathbb{F}_{t}.$$ Thus we can decompose
every element $x$ of $\mathbb{F}_{q}$ $$x=y+\alpha z.$$ The
conjugaison is given by $x\longmapsto x^t$ and transform $\alpha$ in
$\overline{\alpha}$. Therefore $x^t=y+\overline{\alpha}z$, and therefore,
$$x^{t+1}=(y+\alpha z)(y+\overline{\alpha}z)=y^2+(\alpha
+\overline{\alpha})yz+\alpha \overline{\alpha}z^2.$$ From a result
of R. C. Bose and I. M. Chakravarti [1, p.1163 ], we know that the sum
$\alpha +\overline{\alpha}$ as the product $\alpha
\overline{\alpha}$ belong together in $\mathbb{F}_{t}$. Therefore
the form $\phi$ is now a form of degree 2 in $2(n+1)$ variables over
the subfield $\mathbb{F}_{t}$. And its new equation becomes:
$$\phi(y_0,z_0,...,y_n,z_n)= y_0^2+(\alpha
+\overline{\alpha})y_0z_0+\alpha \overline{\alpha}z_0^2+...+
y_n^2+(\alpha +\overline{\alpha})y_nz_n+\alpha
\overline{\alpha}z_n^2.$$
\subsection{Structure of the weights of the codewords}
Let us recall an important propriety on a codes which can be found in [2, pp.10-11].
It is the propriety of divisor of a code. Divisible codes, are interesting because many optimal codes
exhibit  nontrivial divisibility.
\begin{definition} [2, pp.10-11]
We say that a code $\mathcal{C}$ (over any field) is divisible provided all codewords have
weights divisible by an integer $\Delta > 1$. The code is said divisible by $\Delta$; $\Delta$
is called a divisor of the code $\mathcal{C}$, and the largest such divisor is called the divisor
of the code $\mathcal{C}$.
\end{definition}
To our knowledge, during the past forty years since the discovering of the notion of divisor of a code,
the main achievement in the determination of a divisor of a code, has been done for two kind of codes:
cyclic codes over prime fields by R. J. McEliece [16], Griesmer codes (codes meeting the Griesmer bound)
  over prime fields in the binary case by S. M. Dodunekov and N. L. Manev [4]. Making used of the
  divisibility criteria [21, p.323] H. N. Ward [22, p.80, pp.84-87] extend the result to Griesmer codes
  over prime fields in the nonbinary case.
 This is the result on divisibility for Griesmer code.
\begin{theorem} [2, p.86]
Let $\mathcal{C}$ be a linear code over $\mathbb{F}_{p}$, where $p$ is a prime, which meets the Griesmer bound.
Assume that $p^{i}|d(\mathcal{C})$ where $d(\mathcal{C})$ is the minimum distance of $\mathcal{C}$,
then $p^{i}$ is a divisor of the code $\mathcal{C}$.
\end{theorem}
For further details in the study of the divisibility properties of codes, the survey paper of H. N. Ward [23]
 where he generalized the above theorem for Griesmer codes over the field $\mathbb{F}_q$ ($q=p^a$, $p$ is a prime,
 $a$ an integer) by a conjecture [23, p.271] and the recent Ph. D. Thesis of X. Liu [15] where he gave bounds
 on dimension of divisible codes, can be excellent companions.
\begin{theorem} \label{conditiondivisibility}
Let $n$ and $h$ be two positive integers such that $h\leq n$ and $n=sh+r$ where $0\leq r\leq h-1$.
Let us consider the code $C_h(\mathcal{X})$ defined on the non-singular Hermitian variety
$\mathcal{X}:x_0^{t+1}+x_1^{t+1}+...+x_{n-1}^{t+1}+x_n^{t+1}=0$  over the field $\mathbb{F}_q$ ($q=t^2$ and $t=p^a$).
Then  $\Delta = t^{\lambda}$ is the divisor of the code $C_h(\mathcal{X})$ where:
\begin{equation*}
       \lambda =  \begin{cases}
                      n-2       &  \text{if} \qquad  h=2\\

                   2E(\frac{n}{h})-2  &    \text{if}\qquad h \geq 3\qquad \text{and}\qquad r=0 \\

                   2E(\frac{n}{h})-1  &    \text{if}\qquad h \geq 3\qquad \text{and}\qquad h=2r \\
                   2E(\frac{n}{h})+E(\frac{2r}{h})-1 & \text{if}\qquad h \geq 3\qquad \text{and}\qquad h\neq2r
            \end{cases}
\end{equation*}
with E(x) equal the integer part of x.
\end{theorem}
\textbf{Proof:}
Let $f$ be a form of degree $h$ in $n+1$ variables over the field
$\mathbb{F}_{q}$. By using the above transform on the coefficients
of $f$ and the $n+1$ variables, $f$ can newly be written  as:
$$f(x_0,x_1,...,x_{n-1},x_n)= f_0(y_0,z_0,...,y_n,z_n)+
\sum_{i=1}^{h} {\alpha}^i f_i(y_0,z_0,...,y_n,z_n)$$

 where $f_0, f_1,..., f_h$ are
homogeneous polynomials of degree $h$ in $2(n+1)$ variables over the
field $\mathbb{F}_{t}$. By using again the above transform on the
$h-1$ elements ${\alpha}^2$,...,${\alpha}^h$ of $\mathbb{F}_{q}$, we
deduce that finally $$f(x_0,x_1,...,x_{n-1},x_n)=\tilde{f_0}(y_0,z_0,...,y_n,z_n)+ \alpha
\tilde{f_1}(y_0,z_0,...,y_n,z_n)$$  where $\tilde{f_0}$ and
$\tilde{f_1}$ are two homogeneous polynomials of degree $h$ in
$2(n+1)$ variables over the field $\mathbb{F}_{t}$. Therefore the
variety $\mathcal{V}$ defined by one equation on $\mathbb{F}_{q}$,
is defined by two equations on $\mathbb{F}_{t}$ and the variety
 $X\cap\mathcal{V}$ is defined by the following system of three
 equations:

 \begin{equation*}
         \begin{cases}
           y_0^2+(\alpha +\overline{\alpha})y_0z_0+\alpha \overline{\alpha}z_0^2+...+ y_n^2+
           (\alpha +\overline{\alpha})y_nz_n+\alpha \overline{\alpha}z_n^2=0 &    \text{} \\
          \qquad \qquad \qquad \tilde{f_0}(y_0, z_0,...,y_n, z_n)=0 &   \text{} \\
          \qquad \qquad \qquad \tilde{f_1}(y_0, z_0,...,y_n, z_n)=0 & \text
          {}
            \end{cases}
\end{equation*}
which are equations on $\mathbb{F}_{t}$. By the theorem of Ax-Katz,
the number of  common zeros $N$ of $\tilde{f_0}(y_0, z_0,...,y_n,
z_n)$, $\tilde{f_1}(y_0, z_0,...,y_n, z_n)$ and  $\phi$ (in
$\mathbb{F}_{t}^{2(n+1)}$) is divisible by
 $t^{\lambda(S,T,f)}$ where $\lambda(S,T,f)$  is the least non-negative integer such that

\begin{equation}
    \label{Katzinequality}
    \lambda(S,T,f)  \geq  \frac{2(n+1)-(2h+2)}{h}=2(\frac{n}{h}-1).
\end{equation}

 If $h=2$, then  $\lambda(S,T,f)= 2$ and $N$ is divisible by $t^{n-2}$.\\
 Let us suppose now $h\geq 3$ and $n=ah+ r$ where $0\leq r\leq h-1$.\\
 If $r=0$,  then $N$ is divisible by $$t^{2E(\frac{n}{h})-2}.$$
If $r\neq 0$, then $N$ is divisible by

\begin{equation*}
         \begin{cases}
                      t^{2E(\frac{n}{h})-1}  &    \text{if}\qquad h=2r \\
                      t^{2E(\frac{n}{h})+E(\frac{2r}{h})-1} & \text{if} \qquad h\neq2r
            \end{cases}
\end{equation*}
In fact, we have $ 1 <r \leq h-1 $ therefore $ 0 < \frac{2r}{h} < 2 $. When $h=2r$, the least non-negative
integer such that (\ref{Katzinequality}) is verified is $2E(\frac{n}{h})+E(\frac{2r}{h})-2$.  When $h\neq2r$,
we have $\frac{2r}{h} \in ]0,1[ \cup ]1, 2[ $, therefore the least non-negative integer satisfying
(\ref{Katzinequality}) is $2E(\frac{n}{h})+E(\frac{2r}{h})+1-2$.

On the other hands, $\tilde{f_0}(y_0, z_0,...,y_n, z_n)$,
$\tilde{f_1}(y_0, z_0,...,y_n, z_n)$ and  $\phi$
are homogeneous polynomials, therefore $N-1$ is divisible by $t-1$.\\
Let $\mathcal{X}$, $\mathcal{V}_1$ and $\mathcal{V}_2$ be the projective
varieties associated respectively to the forms $\phi$,
$\tilde{f_0}(y_0, z_0,...,y_n, z_n)$, and $\tilde{f_1}(y_0,
z_0,...,y_n, z_n)$, one has $$\#(X \cap \mathcal{V})=\#(X \cap
\mathcal{V}_0\cap\mathcal{V}_1 )= \frac{N-1}{t-1}.$$ Let
$M=\frac{N-1}{t-1}$, and $\lambda$ as above, one has
 $$ M = \frac{jt^{\lambda}-1}{t-1}=j^{\prime}t^{\lambda}+\pi_{{\lambda}-1} \qquad ( \star),$$
  where $j$ and $j^{\prime}$
 are non-null integers  such that $j=j^{\prime}(t-1)+1$.\\
 By the theorem of Ax-Katz again, we get that the number of zeros of the polynomial
 $\phi$ (in  $\mathbb{F}_{t}^{2(n+1)}$) is divisible par $t^n$ and therefore the number
 of zeros of $\mathcal{X}$ in $\mathbb{P}^{2n+1}(\mathbb{F}_t)$ is:
$$ \# \mathcal{X} = \frac{kt^n-1}{t-1}=k^{\prime}t^n+\pi_{n-1}\qquad (\star \star),$$ where $k$ and $k^{\prime}$
 are non-null integers  such that $k=k^{\prime}(t-1)+1$.\\
The weight of the codeword associated to the projective variety $\mathcal{V}$ defined by the form
 $f$ is equal to:
$$w=\# \mathcal{X}-\# (\mathcal{V}\cap \mathcal{X})=\# \mathcal{X}- M, ( \star \star \star).$$
Therefore, from $(\star)$, $(\star \star)$, and $(\star \star \star)$  we deduce that
$$w=k^{\prime}t^n+t^{n-1}+\cdots+t^{{\lambda}+1}+(1-j^{\prime})t^{\lambda}.$$
Thus, $$w \equiv 0 (\mathrm{mod.}\ t^{\lambda} )  .$$
\begin{Remark}
What is important in this paper is that our technique gives directly the divisor,
of the code $\mathcal{X}$, without any knowledge of the minimum distance. In fact most of the
results in the literature relate the determination of a divisor of a code to its minimum distance,
 in the particular cases of cyclic codes, Griesmer codes etc...
Here we don't need to know the minimum distance. We don't need to know if the code is cyclic or attains
the Griesmer bound, but we have a strong information (its divisor).
\end{Remark}
\begin{Remark}
In the case where $h=2$, the first five weights of the code $C_2(\mathcal{X})$ defined on the
non-singular Hermitian variety have been determined by A. Hallez and L. Storme [9, p.9]
in their two tables under the restrictive condition $n< O(t^2)$. There are
divisible by $t^{n-2}$.
\end{Remark}
In the particular case where $h\leq t$, by using Theorem \ref{cardinal hermitienne} and Theorem \ref{TSS}
we deduce that the evaluation map $c$ defining the code $C_h(\mathcal{X})$ is injective. The length and
the dimension of the code $C_h(\mathcal{X})$ on the non-singular Hermitian variety are given respectively by
Theorem \ref{cardinal hermitienne} and relation (\ref{dimducode}). In the general case, a lower bound for
the minimum distance of $C_h(\mathcal{X})$ on the non-singular Hermitian variety has been given by
F. Rodier [17, pp.207-208]. The result of Theorem \ref{conditiondivisibility} gives an improvement of this
lower bound.
\section{The code $C_2(\mathcal{X})$ defined on the Hermitian surface}
The code $C_2(\mathcal{X})$ on the non-singular Hermitian surface
$\mathcal{X}:x_0^{t+1}+x_1^{t+1}+x_{2}^{t+1}+x_3^{t+1}=0$
 has been studied by the first author in
 [5, pp.14-18, pp.28-58], [6], [7]. He found the three first weights of the codewords and their frequency.
  Based on the result obtained in the intersection of quadric surfaces and the non-singular Hermitian surface, he
formulated a conjecture on the fourth and fifth weight of this code.
The resolution of this conjecture depended mainly on how we could improve the upper bound $$2t^3+2t+2$$
found in [7, pp.107-108] for the number of intersection points in the section of the
Hermitian surface and the elliptic quadric.\\
Now we will use a  very simple technique to study it. From the works of R. C. Bose and I. M.
 Chakravarti [1, p.1179], we know that there are exactly $\alpha_t=(t^3+1)(t+1)$ lines
 contained in the non-singular Hermitian surface
 $\mathcal{X}:x_0^{t+1}+x_1^{t+1}+x_{2}^{t+1}+x_3^{t+1}=0$
 and  through each point of $\mathcal{X}$ it pass exactly $t+1$ lines contained in $\mathcal{X}$. We also know
 from [10, p.123] that there is no line in the elliptic quadric surface. Therefore, every line of $\mathcal{X}$
 intersects the elliptic quadric $\mathcal{E}_3$ in at most two points. Thus we deduce that,
 \begin{equation}\label{EllipticHermitian}
\# (\mathcal{X} \cap \mathcal{E}_3)\leq\frac{ 2\alpha_t}{t+1}=2(t^3+1).
\end{equation}
\begin{corollary}\label{surfacedivisibility}
Let $C_2(\mathcal{X})$ be the functional code defined on the non-singular Hermitian surface
$\mathcal{X}:x_0^{t+1}+x_1^{t+1}+x_{2}^{t+1}+x_3^{t+1}=0$ over the field $\mathbb{F}_q$ ($q=t^2$ and $t=p^a$).
Then  $\Delta = t$ is the divisor of  $C_2(\mathcal{X})$.
\end{corollary}
\textbf{Proof:} It is a direct consequence of Theorem \ref{conditiondivisibility}.\\\\
From Corollary \ref{surfacedivisibility} and the relation (\ref{EllipticHermitian}), we deduce a new upper
bound on the number of points in the intersection of the non-singular Hermitian surface and the elliptic quadric:
\begin{equation}\label{EllipticHermitianbest}
\# (\mathcal{X} \cap \mathcal{E}_3)\leq s_4(t)=2t^3+1.
\end{equation}

\begin{theorem}
The fourth weight is $w_4=t^5-t^3+t^2$ and for $t\neq 2$ the corresponding codewords are given by quadrics
which are union of two non-tangent planes to $\mathcal{X}$ and the line of intersection of the two planes
meets $\mathcal{X}$ at a single point. There are exactly $\frac{1}{2}(t-1)t^3(t^4-1)^2$ codewords of fourth
weight.\\
The fifth weight is $w_5=t^5-t^3+t^2+t$ and for $t\neq 2, 3$ the corresponding codewords are given by quadrics
which are union of two non-tangent planes to $\mathcal{X}$ such that the line of intersection of the two
planes  intersects $\mathcal{X}$ in $t+1$ points. There are exactly $\frac{1}{2}(t-1)(t^3+1)(t^2+1)^2t^6$
codewords of fifth weight.
\end{theorem}
\textbf{Proof:}
From the results of tables 1,2, 3  of paragraph 4.2 in [7, pp.111-112], and the improved upper
bound (\ref{EllipticHermitianbest})
obtained on the number of points in the intersection of the Hermitian surface and the elliptic
quadric, we deduce that $s_4(t)$ gives the fourth weight $w_4=t^5-t^3+t^2$.
From the results of tables 1, 2, 3 in [7, pp.111-112] and Corollary \ref{surfacedivisibility},
we deduce that $s_5(t)=2t^3-t+1$ gives the fifth weight $w_5=t^5-t^3+t^2+t$.\\
We will compute now the number of codewords of fourth and fifth weight.
From the fundamental formula of Wan and Yang \lbrack 20\rbrack\ or \lbrack 12, Th. 23.4.3, pp.70-71\rbrack\ ,
we deduce that there are exactly
\begin{equation}\label{Wan2}
N(\mathcal{L}; \mathcal{X}, \Pi_0\mathcal{U}_0)=\frac{t(t^3+1)(t^4-1)}{(t+1)}
\end{equation} lines $\mathcal{L}$ intersecting $\mathcal{X}$ at a single point $\Pi_0$
(i.e. a singular Hermitian variety
 $\Pi_0\mathcal{U}_0$ of rank 1 in $\mathbb{P}^1(\mathbb{F}_q)$). We also know from [1. p.1179] that through the point
 $\Pi_0$, there pass exactly $t+1$ lines which constitute the intersection with $\mathcal{X}$ of the tangent plane at
 $\Pi_0$. Thus, among the $q+1$ planes through the line $\mathcal{L}$, we deduce that there is only one plane which is
tangent to $\mathcal{X}$. Therefore, through the $q(q+1)/2$ pairs of planes through $\mathcal{L}$ we have exactly
\begin{equation}\label{q+2}
   \frac{q(q-1)}{2}
\end{equation}
pairs of planes non-tangent to $\mathcal{X}$.
From (\ref{Wan2})and (\ref{q+2}) we deduce that there are exactly
$(t^2-1).N(\mathcal{L}; \mathcal{X}, \Pi_0\mathcal{U}_0).\frac{q(q-1)}{2}$ codewords of fourth
weight.\\
From \lbrack 20\rbrack\ or \lbrack 12, Th. 23.4.3, pp.70-71\rbrack\ ,  we deduce that there are exactly
\begin{equation}\label{Wan3}
N(\mathcal{L}; \mathcal{X}, \mathcal{U}_1)=\frac{t^4(t^3+1)(t^2+1)}{t+1}
\end{equation} lines $\mathcal{L}$ intersecting $\mathcal{X}$ at $t+1$ points
 (i.e. a non-singular Hermitian variety $\mathcal{U}_1$ in $\mathbb{P}^1(\mathbb{F}_q)$ ).
 Here the $q+1$ planes through the line $\mathcal{L}$ are all
non-tangent to $\mathcal{X}$, because any line of a tangent plane to $\mathcal{X}$ can not
meet $\mathcal{X}$ at $t+1$ points. Therefore, through the $q(q+1)/2$ pair of planes through $\mathcal{L}$
we have exactly
\begin{equation}\label{q+3}
   \frac{q(q+1)}{2}
\end{equation}
pair of planes non-tangent to $\mathcal{X}$.
From (\ref{Wan3}) and (\ref{q+3}) we deduce that there are exactly
$(t^2-1).N(\mathcal{L}; \mathcal{X}, \mathcal{U}_1).\frac{q(q+1)}{2}$
codewords of fifth weight.
\begin{Remark}
If $w_{i}$ $(1\le i\le 5)$ are the first five weights of the code $C_2(\mathcal{X})$,
then there exist degenerate quadrics $\mathcal{Q}$ reaching the Tsfasman-Serre-S\o rensen's
upper bound for hypersurfaces (i.e. $\mathcal{Q}$ is a union of two distinct planes
 $\mathcal{Q}=H_1\cup H_2$), giving codewords of weight $w_i$. And for $i>5$, there is no
quadric which is a union of distinct planes, giving codewords of weight $w_{i}$.
\end{Remark}
\section{Conjecture on the first $2h+1$  weights of the code $C_h(\mathcal{X})$ in
 $\mathbb{P}^{n}(\mathbb{F}_q)$}
The author has also tried to generalize the study of the code $C_h(\mathcal{X})$ to
 $\mathcal{X}$ the non-degenerate Hermitian variety defined by
 $x_{0}^{t+1}+x_{1}^{t+1}+...+x_{n-1}^{t+1}+x_{n}^{t+1}= 0$
  in $\mathbb{P}^{n}(\mathbb{F}_q)$ ($q=t^2$) and conjecture that:\\\\
\textbf{Conjecture } \textit{
\begin{itemize}
\item[{1}] If $w_{i}$ $(1\le i\le 2h+1)$ are the first 2h+1 weights of the code $C_2(\mathcal{X})$,
 then there exist degenerate hypersurface $\mathcal{V}$ reaching the Tsfasman-Serre-S\o rensen's upper bound
  for hypersurfaces (i.e. $\mathcal{V}$ is a union of $h$ distinct hyperplanes  $\mathcal{V}=H_1 \cup...\cup H_h$,
  meeting in a common linear space of codimension 2), giving codewords of weight $w_i$.
\item[{2}]The minimum weight (i.e. $w_1$) codewords are only given by degenerate hypersurfaces which are
 union of $h$ distinct hyperplanes ($\mathcal{V}=H_1\cup...\cup H_h$) such that
 $(H_1\cap...\cap H_h)\cap \mathcal{X}$ is a non-singular Hermitian variety in
 $\mathbb{P}^{N-2}(\mathbb{F}_q)$ and: \\
 \textendash If n is even, the $h$ hyperplanes $H_1$ ,.., $H_h$ are non-tangent to $\mathcal{X}$\\
 \textendash If n is odd, the $h$ hyperplanes $H_1$ , ..., $H_h$ are tangent to $\mathcal{X}$.
\item[{3}] For $i>2h+1$, there is no hypersurface which is a union of distinct hyperplanes, giving
codewords of weight $w_{i}$.
\end{itemize}}
Unfortunately no proof has been found yet.
\begin{Remark}
The conjecture is true for $n=3$ and $h=2$ (see section 4 of this paper and  \lbrack 7,
$\S$ 4.1-4.2\rbrack). The conjecture is also true in the case $n=4$ and $h=2$:
Theorem 4.4 [8, p.143] gives the result. A. Hallez and L. Storme [9] have proved this conjecture
in the case $h=2$ under the condition that $n<O(t^2)$.
 \end{Remark}
\textbf{References}\\
{\footnotesize \lbrack 1\rbrack  \ R. C. Bose and I. M. Chakravarti, Hermitian varieties in finite projective space
$PG(N,q)$. Canadian J. of Math.18 (1966), 1161-1182.\\
\lbrack 2\rbrack  \ W. Cary and V. Pless, Fundamentals of error-correcting codes,
Cambridge University Press, Cambridge 2003.\\
\lbrack 3\rbrack \  I. M. Chakravarti, Some properties and applications of Hermitian varieties in
 finite projective space PG(N,$q^2$) in the construction of strongly regular graphs
 (two-class association schemes) and block designs, Journal of Comb. Theory, Series B,
  11(3) (1971), 268-283.\\
\lbrack 4\rbrack \  S. M. Dodunekov and N. L. Manev, Minimum possible block length of a linear code
for some distance, Problem Inform. Transmission 20, 8-14, (1984).\\
\lbrack 5\rbrack  \ F. A. B. Edoukou, Codes correcteurs d'erreurs construits \`a partir des vari\'et\'es
alg\'ebriques. Ph.D Thesis, Universit\'e de la M\'editerran\'ee (Aix-Marseille II), France, (2007).\\
\lbrack 6\rbrack  \ F. A. B. Edoukou, Codes defined by forms of degree 2 on Hermitian surface
 and S\o rensen conjecture. Finite Fields and Their Applications 13 (3), 616-627, (2007).\\
\lbrack 7\rbrack  \ F. A. B. Edoukou, The weight distribution of the functional codes defined by
 forms of degree 2 on Hermitian surfaces. Journal de Th\'eorie de Nombres de Bordeaux 21 (1),
 105-107, (2009).\\
\lbrack 8\rbrack  \ F. A. B. Edoukou, Codes defined by forms of degree 2 on non-degenerate Hermitian
varieties in $\mathbb{P}^{4}(\mathbb{F}_q)$. Designs Codes and Cryptography 50, 135-146, (2009).\\
\lbrack 9\rbrack  \ A. Hallez and L. Storme, Functional codes arising from quadric intersections
with Hermitian varieties. 10 pages, Submitted to Finite Fields and Theirs Applications, (2009).\\
\lbrack 10\rbrack  \ J. W. P. Hirschfeld, Projective Geometries Over Finite Fields (Second Edition)
 Clarendon  Press. Oxford 1998.\\
\lbrack 11\rbrack  \ J. W. P. Hirschfeld, Finite projective spaces of three dimensions,
 Clarendon press. Oxford 1985. \\
\lbrack 12\rbrack  \ J. W. P. Hirschfeld,  General Galois Geometries, Clarendon press. Oxford 1991.\\
\lbrack 13\rbrack  \ N. M. Katz,  On a Theorem of Ax, American J. of Mathematics 93 (1971),
 no 2, 485-499.\\
\lbrack 14\rbrack  \ G. Lachaud, Number of points of plane sections and linear codes defined on
 algebraic varieties;  in " Arithmetic, Geometry, and Coding Theory ". (Luminy, France, 1993),
 Walter de Gruyter, Berlin-New York, (1996), 77-104. \\
\lbrack 15\rbrack  \ X. Liu, On divisible codes over finite fields, Ph. D. Thesis
California Institute of Technology, Pasadena, California, USA, (2006).\\
\lbrack 16\rbrack \  R. J. McEliece, Weight congruences for p-ary codes, Discrete Math.
 3, 177-192, (1972).\\
\lbrack 17\rbrack \  F. Rodier, Codes from flag varieties over a finite field, Journal of Pure
and Applied Algebra 178 (2003), 203-214.\\
\lbrack 18\rbrack  \ J. -P. Serre, Lettre \`a M. Tsfasman, In "Journ\'ees Arithm\'etiques de
Luminy (1989)", Ast\'erisque 198-199-200 (1991), 3511-353.\\
\lbrack 19\rbrack  \ A. B. S\o rensen, Rational points on hypersurfaces, Reed-Muller codes and
 algebraic-geometric codes. Ph. D. Thesis, Aarhus, Denmark, 1991.\\
\lbrack 20\rbrack \ Z.-X. Wan and B. -F. Yang, Studies in finite geometries and the construction of
incomplete block designs. III. Some 'Anzahl' theorems in unitary geometry over finite fields and
their applications. Chinese Math.7, (1965), 252-264.\\
\lbrack 21\rbrack \  H. N. Ward, Weight polarization and divisibility, Discrete Mathematics 83, 315-326,
(1990).\\
\lbrack 22\rbrack \  H. N. Ward, Divisibility of codes meeting the Griesmer Bound, Journal of Combinatorial
Theory, Series A 83, 79-93 (1998).\\
\lbrack 23\rbrack \  H. N. Ward, Divisible codes - a survey, Serdica Mathematical
 Journal 27, 263-278 (2001).}
\end{document}